\documentclass[a4paper,10pt]{article}
\usepackage{amsmath,amssymb, amsthm,epsfig,xypic}
\usepackage[all]{xy}

\input xy
\xyoption{all}

\newtheorem{theo}{\bf{Theorem}}[section]
\newtheorem{lem}[theo]{Lemma}

\newtheorem{prop}[theo]{Proposition}
\newtheorem{defi}[theo]{Definition}

\begin{document}
\title{Counting overlattices in automorphism groups of trees}
\author{Seonhee Lim}
\date{}
\maketitle
\begin{abstract}
We give an upper bound for the number $u_{\Gamma}(n)$ of ``overlattices'' in the automorphism group of a tree, containing
 a fixed lattice $\Gamma$ with index $n$. For an example of $\Gamma$ in the automorphism group of a $2p$-regular tree
  whose quotient is a loop, we obtain a lower bound of the asymptotic behavior as well.
\\

\medskip

Nous donnons une borne sup\'erieure pour le nombre $u_{\Gamma}(n)$ de ``surr\'eseaux'' contenant un r\'eseau fix\'e
d'indice $n$ dans le groupe d'automorphismes d'un arbre. Dans le cas d'un arbre $2p$-r\'egulier $T$, et d'un r\'eseau
$\Gamma$ tel que $\Gamma \backslash T$ soit une boucle, nous obtenons aussi une minoration du comportement asymptotique.

\end{abstract}

\paragraph{Introduction.}
Given a connected semisimple Lie group $G$, the Kazhdan-Margulis lemma says that there exists a positive lower bound for
the covolume of cocompact lattices in $G$. This is no longer true when $G$ is the automorphism group of a locally finite
tree. Bass and Kulkarni (for cocompact lattices, see \cite{BK}) and Carbone and Rosenberg (for arbitrary lattices in
uniform trees, see \cite{CR}) even constructed examples of increasing sequences of lattices
$(\Gamma_i)_{i \in \mathbb{N}}$ in $Aut(T)$ whose covolumes tend to 0 as $i$ tends to $\infty$.

If $\Gamma$ is a cocompact lattice in the group $Aut(T)$ of automorphisms of a locally finite tree $T$, there is only a
finite number $u_\Gamma (n)$ of ``overlattices'' $\Gamma'$ containing $\Gamma$ with fixed index $n$ (\cite{B}). Thus a
 natural question, which was raised by Bass and Lubotzky (see \cite{BL}), would be to find the asymptotic behavior of
 $u_\Gamma(n)$ as $n$ tends to $\infty$.

In \cite{G}, Goldschmidt proved that there are only 15 isomorphism classes of (3,3)-amalgams. Thus for lattices $\Gamma$
in the automorphism group of a 3-regular tree $T$ whose edge-indexed quotient is \;\;
\begin{picture}(0,0)%
\includegraphics{essai3.pstex}%
\end{picture}%
\setlength{\unitlength}{1184sp}%
\begingroup\makeatletter\ifx\SetFigFont\undefined%
\gdef\SetFigFont#1#2#3#4#5{%
  \reset@font\fontsize{#1}{#2pt}%
  \fontfamily{#3}\fontseries{#4}\fontshape{#5}%
  \selectfont}%
\fi\endgroup%
\begin{picture}(2341,623)(968,-444)
\put(1276,-136){\makebox(0,0)[lb]{\smash{\SetFigFont{7}{8.4}{\rmdefault}{\mddefault}{\updefault}{$3$}%
}}}
\put(2701,-136){\makebox(0,0)[lb]{\smash{\SetFigFont{7}{8.4}{\rmdefault}{\mddefault}{\updefault}{$3$}%
}}}
\end{picture}
\;\;,
\;\;one has $u_\Gamma(n)=0$ for $n$ big enough. Moreover it is conjectured by Goldschmidt and Sims that there is only a
 finite number of (isomorphism classes of) $(p,q)$-almalgams, for any prime numbers $p$ and $q$.

\paragraph{}
 In this paper, we give two results: an upper bound of $u_\Gamma(n)$ for any cocompact lattice, and a surprisingly big
exact asymptotic growth of $u_\Gamma(n)$ for a specific lattice $\Gamma$ in the automorphism group of a
$2p$-regular tree.

\begin{theo} Let $\Gamma$ be a cocompact lattice in $Aut(T)$. Then there are some positive constants $C_0$ and $C_1$
depending on $\Gamma$, such that
\begin{equation*}
\forall n \geq 1, \;\;\;\;\;\;\; u_{\Gamma} (n) \leq C_0 n^{C_1 \log^2(n)}.
\end{equation*}
\end{theo}

\begin{theo} Let $p$ be a prime number and let $T$ be a $2p$-regular tree. Let $\Gamma$ be a cocompact lattice in $Aut(T)$
such that the quotient graph of groups is a loop whose edge stabilizer is trivial and whose vertex stabilizer is
a finite group of order $p$.
$$
\centerline{\epsfbox{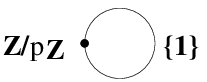}}
$$
Let $n=p_0 ^{k_0} p_1 ^{k_1} \cdots p_t ^{k_t}$ be the prime decomposition of $n$ with $p_0=p$. Then there exist positive
 constants $c_0, c_1$ such that $\underset{k_0 \to \infty}{\limsup} \frac{u_\Gamma(n)}{n^{c_1 \log n}} \leq c_0$.
For $n=p_0^{k_0} (k_0 \geq 3)$, we also have $u_\Gamma(n)\geq n^{\frac{1}{2}(k_0-3)}$.
\end{theo}

It is easy to see (\cite{B}) that $$u_{\Gamma}(n) \leq \sum_{{\underset{\Gamma'' \subset
\Gamma}{[\Gamma:\Gamma''] | n!}}} |N_{Aut(T)}(\Gamma '')/ \Gamma''|,$$ thus we could hope to use the results of
Lubotzky on subgroup growth (see for instance \cite{L1}, \cite{L2}). However, the estimations given in this way
do not seem to be sharp enough. Thus our strategy consists in using the correspondence between cocompact
lattices and graphs of groups (the Bass-Serre theory, see section 1) and reducing the problem to counting
certain isomorphism classes of covering graphs of groups of index $n$ (see section 2).

Together with the sharply contrasting examples satisfying the Goldschmidt-Sims conjecture, the examples in Theorem 0.2
are presently the only known behavior for overlattice counting functions.

\medskip

\textit{Acknowledgements:} We thank Alex Lubotzky for introducing the subject and the problem as well as for helpful
discussions. We thank Gregory Margulis for his guidance and L\'aszl\'o Pyber and Gabe Rosenberg for explaining their
works (\cite{P}, \cite {CR}) and \cite {BK}. Finally, we are grateful to Fr\'ed\'eric Paulin for his constant help and
encouragement.

\section{Preliminaries}
In this section, we briefly recall some background on group actions on trees and the theory of graphs of groups,
and we explain the correspondence between overlattices and coverings of graphs of groups. We refer the reader to
\cite{S}, \cite{B} and \cite{BL} for details on the standard material, gathered in section 1.1.

\paragraph{}
Throughout the paper, we denote by $T$ a locally finite tree, i.e., a tree having finite valence at each vertex. We
denote by $Aut(T)$ the group of automorphisms without inversions of the tree $T$.
A subgroup $\Gamma$ of $Aut(T)$ is \textit{discrete} if the stabilizer $\Gamma_{x}$ is finite for some, thus for every,
vertex $x$ of $T$. The {\it covolume} of $\Gamma$ is defined by
$$Vol\;(\Gamma \backslash \backslash T)= \underset{x \in \Gamma \backslash VT}{\sum} \frac {1}{|\Gamma_{x}|}.$$

A discrete subgroup is a \textit{lattice} if its covolume is finite. In this case, $Aut(T)$ is unimodular, and
the covolume is equal (up to a constant depending only on $T$) to the volume of $\Gamma \backslash Aut(T)$
induced by the Haar measure on the locally compact group $Aut(T)$ \cite{BL}.
 A lattice $\Gamma$ is ca
 lled {\it cocompact} if the quotient graph $\Gamma \backslash T$ is finite.
 An \textit{overlattice} of $\Gamma$ is a lattice of $Aut(T)$ containing $\Gamma$ with finite index.

\paragraph{1.1. Cocompact lattices and finite graphs of finite groups}

By a graph of groups $(X,G_\bullet )$, we mean a connected graph $X$, groups $G_{x}$ and
$G_{e}=G_{\overline{e}}$ assigned to each vertex $x$ in $VX$ and each edge $e$ in $EX$, together with injections
$G_e \rightarrow G_x$ for each edge $e$ with origin $o(e)=x$. This injective map will be denoted by $\alpha_e$,
whatever the graph of group is. The {\it edge-indexed graph} of the graph of groups $(X, G_\bullet)$ is the
graph $X$ with {\it index} $i(e)= |G_{o(e)}|/ |G_e|$ associated to each edge $e$.

\medskip

\paragraph{} To every subgroup $\Gamma$ of $Aut(T)$ is associated a graph of groups, well-defined up to isomorphism of
graph of groups (see definition below), whose graph is the quotient graph $\Gamma \backslash T$. We will call it
a quotient graph of groups of $\Gamma$ and denote it by $\Gamma \backslash\backslash T$. According to
\cite{B}(section 3), a construction of $\Gamma \backslash\backslash T$ proceeds as follows.
 Let $p: T \to X$ be
the canonical projection.

Choose subtrees $R \subset S \subset T$ such that $p |_R :R \to X$
is bijective on vertices, $p|_S :S \to X$ is bijective on edges,
and for each edge $e$ in $E(S)$, at least one of $o(e)$, $t(e)$
belongs to $R$. Define $\tilde{x} =p|_R^{-1} (x)$ for each $x$ in
$VX$ and $\tilde{e} =p|_S^{-1}(e)$ with
$\bar{\tilde{e}}=\tilde{\bar{e}}$ for each $e$ in $EX$. For each
$x$ in $VS$, choose an element $g_x$ in $\Gamma$ such that $g_x
x=\widetilde{p(x)}$. We can choose $g_x=1$ for all $x$ in $VR$.
Now let $G_x$ be the stabilizer $\Gamma_{\tilde{x}}$ of $\tilde x$
in $\Gamma$ for $x$ in $VX \cup EX$. The injective map $\alpha_e$
is defined as
\begin{align*}
\alpha_e : G_e &\; \longrightarrow \; G_{o(e)}\\
            s    &\;  \longmapsto  \; g_{o(\tilde{e})} s g^{-1}_{o(\tilde{e})}.
\end{align*}
 Note that each $\alpha_e$ is merely an inclusion for $e$ such that $o(e)$ is a vertex of $VR$.

Conversely, for any graph of groups $(X, G_\bullet)$, there exists a tree $T$ and a group $\Gamma$ acting on the tree
$T$ (unique up to equivariant tree isomorphism) such that $(X, G_\bullet)$ is isomorphic to
$\Gamma \backslash\backslash T$.
Let us call $(T, \Gamma)$ a \textit{universal cover} of $(X, G_\bullet)$ and $\Gamma$ its {\it fundamental group}.

Fix $x_0 \in VX$. The fundamental group $\Gamma$ of $(X, G_\bullet)$ based at $x_0$ is defined as follows. The
\textit{path group} $\Pi(X, G_\bullet)$ is defined
 by
\begin{equation*}
 \left(\underset{x \in VX}{*}G_{x}\right)*F(EX)/\langle e^{-1}=\overline{e}, e\alpha_{\overline{e}}(g)e^{-1}=
 \alpha_e(g) :g \in G_{e}\rangle ,
\end{equation*}
where $F(EX)$ denotes the free group with basis $EX$. For $x, x'$ in $VX$, we denote by $\pi[x,x']$ the subset
of $\Pi(X,G_\bullet)$ which consists of elements of the form $g_0 e_1 g_1 e_2 \cdots g_{n-1} e_n g_n$ where
$e_i$ is an edge from $x_{i-1}$ to $x_{i}$, $g_i \in G_{x_i}$ and $x_0=x, x_n = x'$. The \textit{fundamental
group of $(X, G_\bullet)$ based at $x_0$} is $\Gamma=\pi_1(X,G_\bullet,x_0)= \pi[x_0, x_0]$, endowed with the
group structure induced by $\Pi(X, G_\bullet)$.

The universal cover
 $\widetilde{(X, G_\bullet, x_0)}$ of $(X,G_\bullet)$ based at $x_0$ is defined as follows. It has as vertex set
\begin{equation*}
V(\widetilde{(X, G_\bullet, x_0)})=\underset{x \in VX}{\coprod} \pi[x_0,x]/G_x,
\end{equation*}
and there is an edge between two distinct points $[g]$ in $\pi[x_0, x]/G_x$ and $[g']$ in $\pi[x_0, x']/G_{x'}$
if and only if $g^{-1}g' \in G_x e G_{x'}$ where $e$ is an edge in $X$ from $x$ to $x'$. The fundamental group
$\pi_1(X, G_\bullet, x_0)=\pi[x_0, x_0]$ acts on $\widetilde{(X, G_\bullet, x_0)}$ by the natural left action.
The graph $\widetilde{(X, G_\bullet, x_0)}$ is a tree and moreover, for any other universal cover $(T, \Gamma)$
of $(X, G_\bullet)$, there is an isomorphism $\psi$ between $\Gamma$ and $\pi_1 (X, G_\bullet, x_0)$ and a
$\psi$-equivariant graph isomorphism between $T$ and $\widetilde{(X, G_\bullet, x_0)}$, see for example
\cite{S}.

\medskip
A graph of groups is called \textit{faithful (or effective)} if there is no edge subgroup family $(N_{e})_{e \in
EX}$
 satisfying the following conditions:
\begin{enumerate}

\item[i)] for each $e$ and $e'$ in $EX$ such that $o(e)=o(e')$, the images of $N_e$ and $N_{e'}$ coincide:
$\alpha_e (N_e)=\alpha_{e'} (N_{e'})$. Let us denote it by $N_{o(e)}$.
\item[ii)] For each $x$ in $VX$, $N_x$ is a nontrivial normal subgroup in $G_x$.
\end{enumerate}

It is shown in \cite{B} that the graph of groups $(X, G_\bullet)$ is faithful if and only if its fundamental
group $\Gamma$ is a subgroup of $Aut(T)$ for $T$ its universal cover, i.e., if and only if the map $\Gamma
\longrightarrow Aut(T)$ is injective. The fundamental group of a faithful finite graph of finite groups is a
cocompact lattice in the automorphism group of its universal covering tree and conversely, a quotient graph of
groups of a cocompact lattice in the automorphism group of a locally finite tree is a faithful finite graph of
finite groups.

\paragraph{}In \cite{B}, Bass defines a covering of graphs of groups in such a way that the induced map between the
corresponding fundamental groups is a group monomorphism.

\begin{defi} Let $(X,G_\bullet)$ and $(Y,H_\bullet)$ be two graphs of groups.
We call a \textit{morphism of graphs of groups}, which we denote
by $\phi_\bullet=(\phi, \phi_x, \gamma_x):(X,G_\bullet)
\rightarrow (Y, H_\bullet)$, the following data

\begin{enumerate}
\item[(i)] a graph morphism $\phi:X \rightarrow Y$,
\item[(ii)] group homomorphisms $\phi_{x}:G_{x} \rightarrow H_{\phi(x)}$ and $\phi_{e}:G_{e} \rightarrow H_{\phi(e)}$,
for every vertex $x$ and every edge $e$ of $X$,
\item[(iii)] families of elements $(\gamma_x) _{x \in VX} \in \pi_1(Y, H_\bullet, \phi(x))$ and $(\gamma_e)_{e \in EX}
\in \Pi(Y,H_{\bullet})$
\end{enumerate}
such that
\begin{enumerate}
\item[] for every edge $e$ of $X$ with origin $x$, we have $\gamma_x^{-1}\gamma_e \in H_{\phi(x)}$
and the following diagram commutes. Here $ad(g)(s)=gsg^{-1}$.
$$\xymatrix{
G_e \ar[d]^{\phi_e} \ar[r]^{\alpha_e} & G_x \ar[d]^{\phi_x}\\
H_{\phi(e)} \ar[r]^-{ad(\gamma_x^{-1} \gamma_e) \circ \alpha_{\phi(e)}} & H_{\phi(x)}}$$

\end{enumerate}
\end{defi}

The \textit{induced homomorphism of path groups}
$\Phi=\Phi_{\phi_{\bullet}}: \Pi(X,G_{\bullet}) \to
\Pi(Y,H_{\bullet})$, is defined as follows on generators (see
\cite{B}): $\Phi(g)=\gamma_x\phi_x(g)\gamma_x^{-1}$ for $g \in
G_x$ and $x \in VX$, $\Phi(e)=\gamma_e\phi(e)\gamma_{\bar e}^{-1}$
for $e \in EX$. The induced homomorphism on path groups restricts
to a homomorphism $\pi_1(X,G_{\bullet},x_0) \to
\pi_1(Y,H_{\bullet},\phi(x_0))$, which we will denote again by
$\Phi$.

The induced homomorphism
$\Phi=\Phi_{\phi_\bullet}:\pi_1(X,G_\bullet,x_0)\rightarrow
\pi_1(Y,H_\bullet,\phi(x_0))$ gives a $\Phi_{x_o}$-equivariant
graph isomorphism
$\tilde{\phi}:\widetilde{(X,G_\bullet,x_o)}\rightarrow
\widetilde{(Y,H_\bullet, \phi(x_0))}$ defined by
\begin{align*}
[g] \in \pi[x_0,x]/G_x \mapsto [\Phi(g)\gamma_x] \in
\pi[\phi(x_0), \phi(x)]/H_{\phi(x)}.
\end{align*}

A morphism $\phi_{\bullet}=(\phi, \phi_x, \gamma_x)_{x \in VX \cup EX}$ of graphs of groups is an {\it
isomorphism of graphs of groups} if $\phi$ is a graph isomorphism and $\phi_x$ are all group isomorphisms. In
this case, $\phi_\bullet ^{-1}=(\phi^{-1}, \phi'_y, \gamma'_y)$ where $\phi'_y=\phi_{\phi^{-1}(y)}$ and
$\gamma'_y=\Phi^{-1}(\gamma_{\phi^{-1}(y)})^{-1}$ for $y \in VY \cup EY$.

\begin{defi} A morphism of graphs of groups $\phi_{\bullet}$ is furthermore called a \textit{covering} if
\begin{enumerate}
\item[(a)] the maps $\phi_e$ and $\phi_x$ are injective for all $x$ and $e$,

\item[(b)]for every edge $f$ of Y with origin $\phi(x)$,where $x$ is in $VX$, the well-defined map
\begin{align*}
 \Phi_{x/f} : \coprod_{e \in \phi^{-1}(f), o(e)=x}
G_x /\alpha_e(G_e) \longrightarrow &\; H_{\phi(x)} /\alpha_f (H_f)\\
      [g]_e \longmapsto &\; [\phi_x (g) \gamma^{-1}_x \gamma_e ]_f
\end{align*}
is bijective.
\end{enumerate}
\end{defi}

\medskip

By the condition $(b)$ in Definition 1.2, we have $\underset{e \in \phi^{-1}(f), o(e)=x}\sum \frac{|G_x|}{|G_e|}
= \frac{|H_{\phi(x)}|}{|H_{f}|}$ for every edge $f$ of $Y$ with origin $\phi(x)$.
 Summing over all vertices $x$ such that $\phi(x)=y$, it follows that the value of $n: =\underset{x \in \phi^{-1}(y)}
 \sum \frac{|H_y|}{|G_x|} = \underset{e \in \phi^{-1}(f)} \sum \frac{|H_f|}{|G_e|}$ does not depend on vertices and
 edges, since the graph $Y$ is connected. Note that $n$ is an integer since $\phi_x(G_x)$ is a subgroup of $H_y$ for
 each $x$ such that $\phi(x)=y$. A covering graph of groups with the above $n$ is said to be {\it $n$-sheeted}.

Note also that by the condition $(b)$, a covering of graphs of groups induces a covering of the corresponding
edge-indexed graphs. Recall that a {\it covering} $\phi: (X,i) \to (Y, i)$ of edge-indexed graphs is a graph
morphism $\phi$ such that $ \sum_{e \in \phi^{-1}(e'), o(e)=x} i(e)=i(e')$.

\begin{theo}[\cite{B}, Prop. 2.7] The morphism $\phi_{\bullet}$ is a covering if and only if
 $\;\Phi: \pi_1(X,G_{\bullet},x_0) \to \pi_1(Y,H_{\bullet},\phi(x_0))$ is injective and $\widetilde \phi :
\widetilde{(X,G_{\bullet},x_0)} \to \widetilde{(Y,H_{\bullet},\phi(x_0))}$ is an isomorphism.\end{theo}

\paragraph{1.2. Counting overlattices}
  Let $\Gamma$ be a cocompact lattice in $Aut(T)$. Set $$U(n)=U_\Gamma (n)=\{\Gamma': \Gamma \subset \Gamma'\subset
Aut(T),\; [\Gamma':\Gamma]=n\}$$ and let $u(n)=u_\Gamma (n)=|U(n)|$ be the number of overlattices of $\Gamma$ of
index $n$. It is shown in \cite{BK} that $u(n)$ is finite. We are interested in the asymptotic behavior of
$u(n)$. For that purpose, we will show in this section that there is a bijection between overlattices of
$\Gamma$ and isomorphisms classes of coverings of graphs of groups by the quotient graph of groups of $\Gamma$,
in the following sense.

\begin{defi}
Let $\phi_\bullet=(\phi,\phi_x,\gamma_x): (X,G_\bullet) \to(Y,H_\bullet)$ and
$\psi_\bullet=(\psi,\psi_x,\gamma'_x): (X,G_\bullet) \to(Y',H'_\bullet)$ be two coverings of graphs of groups.
An isomorphism between them consists of a pair $\{\theta_\bullet=(\theta, \theta_y,\rho_y): (Y,H_\bullet) \to
(Y',H'_\bullet), (\zeta_x)_{x \in VX \cup EX}\}$ where $\theta_\bullet$ is an isomorphism of graphs of groups
$(Y,H_\bullet) \to (Y',H'_\bullet)$ and $(\zeta_x) \in H'_{\psi(x)}$ are such that
\begin{enumerate}
\item[a)] $\theta \circ \phi=\psi $ as a map of graphs,
\item[b)] For any $x \in VX \cup EX$, we have $\psi_x=Ad(\zeta_x)\theta_{\phi(x)} \circ \phi_x$
as maps $G_x \to H'_{\psi(x)}$,
\item[c)] $\gamma'_x=\Theta(\gamma_x)\rho_{\phi(x)}\zeta_x$
for any $x \in VX \cup EX$.
\end{enumerate}
\end{defi}

 For a given overlattice $\Gamma'$ of $\Gamma$, we can construct a covering $m^{\Gamma'}$ of graphs of groups as follows.
 Let $Y = \Gamma' \backslash T$ and $p' : T \to Y$ be the canonical projection.

\medskip

 Define subtrees $R'$ and $S'$ of $R$ and $S$, respectively, in the following
 way. For each vertex $y$ of $Y$, choose one vertex from each set $\{{p'}^{-1}(y)\} \cap VR$ and call it $\tilde{y}$. Let
 $R'$ be the subgraph of $R$ with vertices $\{\tilde{y}: y \in Y\}$. Since $R$ is a tree, we can choose vertices
 $\tilde{y}$ so that $R'$ is connected. Let $S'$ be the maximal subtree of $S$ containing $R'$ such that
 $p'|_{S'}$ is injective on the edges. Choose elements $g'_x \in \Gamma'$ such that $g'_x x =\widetilde{p'(x)}$.
 The graph of groups $(Y, H_\bullet)$ is defined with respect to $R'$, $S'$ and $g'$'s as $(X, G_\bullet)$ is defined in
 section 1.1.

 Now the covering of graphs of groups, which will be denoted by $m=m^{\Gamma'}:(X,G_\bullet) \to (Y, H_\bullet)$ is defined
 as follows. For the graph morphism $m:X \to Y$, take the natural projection $\pi$. For the group morphisms
 $m_x: G_x \to H_{m(x)}$,take an
 element $\sigma_x$ in $\Gamma'$ which sends $\tilde{x}$ to $\widetilde{p(x)}$. We can choose $\sigma_x=1$ if
 $\tilde x \in VR' \cup ES'$. Note that $p(x)$ is a vertex of
 $Y$, thus $\widetilde{p(x)} \in R'$ whereas $x$ is a vertex of $X$, thus $\tilde x \in R$.

   Let $m_x=ad(\sigma_x)\circ \iota$ be
 the injection followed by the conjugation ($g \mapsto \sigma_x g \sigma_x^{-1}$).  Since $G_x$ stabilizes $\tilde{x} \in VT \cup ET$,
  the group $\sigma_x G_x \sigma_x^{-1}$ stabilizes $\widetilde{p(x)} \in VT \cup ET$, thus it is a subgroup of $H_{p(x)} =
 \Gamma'_{\widetilde{p(x)}}$, for $x \in VX \cup EX$. For the elements $\gamma_x, \gamma_e$ in (iii) of
 Definition 1.1, take $\gamma_x=\sigma_x^{-1}$ and $\gamma_e=g_e \sigma_e^{-1} g'^{-1}_{m(e)}$. It follows that
\begin{align*}
ad(\gamma_x^{-1} \gamma_e) \circ \alpha_{m(e)} \circ m_e &= ad(\gamma_x^{-1} \gamma_e) \circ ad(g'_m(e)) \circ
ad(\sigma_e)\\                                           &=ad(\sigma_x g_e \sigma_e^{-1} g'^{-1}_{m(e)}) \circ
ad(g'_{m(e)}) \circ ad(\sigma)\\                          &=ad(\sigma_x g_e) = ad(\sigma_x)\circ ad(g_e)= m_{x}
\circ \alpha_e.
\end{align*}

Since $\gamma_x$'s are the elements of $\Gamma'$, the map $m^{\Gamma'}$ is indeed a covering of graphs of
groups.

\begin{prop}Let $\Gamma$ be a cocompact lattice of $Aut(T)$ and $(X, G_\bullet)$ be its quotient graph of groups. The
map $\Gamma'\mapsto m^{\Gamma'}$ induces a bijection $m$ between the set of overlattices of $\Gamma$ of index
$n$ and the set of isomorphism classes of the $n$-sheeted coverings of faithful graphs of groups by $(X,
G_\bullet)$.
\end{prop}

The following lemma shows that the map $m: \Gamma \mapsto m^\Gamma$ is well-defined.
\begin{lem}\label{L:1} Let $\Gamma$ be a lattice in $T$, and let $\Gamma' \supset \Gamma$ be an overlattice. Fix
$(R,S,g_e)$ giving rise to a graph of groups structure $(X,G_\bullet)$ on $\Gamma \backslash T$. Let
$(R',S',g'_e)$ (resp. $(R'',S'',g''_e)$) be a data giving rise to a graph of groups structure $(Y,H_\bullet)$
(resp. $(Y',H'_\bullet)$) on $\Gamma' \backslash T$, and let $(\theta'_x)_{x \in VX \cup EX}$ (resp.
$(\theta''_x)_{x \in VX \cup EX}$) be a data giving rise to a covering $\phi_\bullet=(\phi,\phi_x,\gamma'_x):
(X,G_\bullet) \to (Y,H_\bullet)$ (resp. $\psi_\bullet=(\psi,\psi_x,\gamma''_x): (X,G_\bullet) \to
(Y',H'_\bullet)$ ). Then the two coverings $\phi_\bullet$ and $\psi_\bullet$ are isomorphic.\end{lem}

\noindent \textit{Proof.} Recall that by definition, we have $\theta'_x: \tilde x \mapsto\widetilde{\phi(x)}$
and $\theta''_x: \tilde x \mapsto \widetilde{\psi(x)}$, where $\theta'_x$ and $\theta''_x$ are in $\Gamma'$.
Recall also that $\gamma'_x={\theta'_x}^{-1}, \gamma''_x={\theta''_x}^{-1}$ for $x \in VX$ and $\gamma'_e=g_e
{\theta'_e}^{-1} {g'_{\phi(e)}}^{-1}$, $\gamma''_e=g_e{\theta''_e}^{-1} {g''_{\psi(e)}}^{-1}$ for $e \in EX$.
Now we want to construct an isomorphism $\{\theta_\bullet:(Y,H_\bullet) \to (Y',H'_\bullet), \zeta_x\}$ of
covering of graph of groups. First notice that there is a canonical bijection $\theta: Y \simeq \Gamma'
\backslash T \simeq Y'$. It lifts to a bijection $\tilde \theta: R' \to R''$ and it extends to a unique
bijection $\tilde \theta: S' \to S''$. Let us choose arbitrary elements $\xi_y \in \Gamma'$ for $y \in VY \cup
EY$ such that $\xi_y(\tilde y) =\widetilde{\theta(y)}$ and define maps
$$\theta_y: H_y =\Gamma'_{\tilde y} \to \Gamma'_{\widetilde{\theta(y)}}=H'_{\theta(y)}, h \mapsto \xi_y h \xi_y^{-1}.$$
We have a morphism of graphs of groups $\theta_\bullet=(\theta, \theta_y, \sigma_y): (Y,H_\bullet) \to
(Y',H'_\bullet)$ by setting $\sigma_y=\xi_y^{-1}$ for $y \in VY$ and $\sigma_e={g'_e}^{-1}{\xi_e}^{-1}
{g''_{\theta(e)}}$ for $e \in EY$. It is clear by construction that this is an isomorphism of graphs of groups
(all maps are isomorphisms of groups). Finally, put $\zeta_x=\xi_x{\theta'_x}{\theta''_x}^{-1}$. We compute, for
$g \in G_e$, $e \in EX$:
$$\psi_e(g)=\theta''_e g {\theta''_e}^{-1},$$
\begin{equation}
\begin{split}
Ad(\zeta_e^{-1}) \theta_{\phi(e)}\phi_e(g)&=\zeta_e^{-1}\xi_{\phi(e)}\theta'_e
g {\theta'_e}^{-1} {\xi_{\phi(e)}}^{-1}\zeta_e\\
&={\theta''_e} g {\theta''_e}^{-1}
\end{split}
\end{equation}
as desired. A similar computation holds for $\psi_x: G_x \to H'_{\psi(x)}$ when $x \in VX$. This proves
condition (b) in the definition of isomorphism of coverings. Condition (c) follows from the very definition of
$\zeta_x, \sigma_y, \gamma'_x$ and $\gamma''_x$. \qed

Let us define the inverse map $\phi_\bullet \mapsto \Gamma_\phi$ of $m$ as follows. Set
$\Gamma_Y:=\pi_1(Y,H_\bullet,\phi(x_0)) \subset Aut((\widetilde{Y,H_{\bullet},\phi(x_0)}))$. We
 define an embedding $i_{\phi} : \Gamma_Y \to Aut((\widetilde{X,G_{\bullet},x_0}))$ as follows :
$$i_{\phi}(u) \cdot v=\tilde{\phi}^{-1}(u \cdot
\tilde{\phi}(v))\qquad for\; u \in \Gamma_Y\; and\; v \in V(\widetilde{X,G_{\bullet},x_0}) \cup
E(\widetilde{X,G_{\bullet},x_0}).$$ Let us denote by $\Gamma_{\phi} \subset
Aut((\widetilde{X,G_{\bullet},x_0}))$ the image of $i_{\phi}$. The following lemma shows that this map is
well-defined.

\begin{lem} If $\phi_\bullet:(X,G_\bullet) \to(Y,H_\bullet)$ and $\psi_\bullet:(X,G_\bullet) \to (Y',H'_\bullet)$ are
isomorphic coverings of graphs of groups, then the corresponding subgroups $\Gamma_\phi \subset
Aut((\widetilde{X,G_\bullet,x_0}))$ and $\Gamma_\psi \subset Aut((\widetilde{X,G_\bullet,x_0}))$
coincide.\end{lem}

\noindent \textit{Proof.} We have a triangle of morphisms of path groups
$$\xymatrix{
\Pi(X,G_\bullet) \ar[r]^-{\Phi} \ar[dr]^-{\Psi}& \Pi(Y,H_\bullet) \ar[d]^-{\Theta}\\
& \Pi(Y',H'_\bullet) }$$ We claim that this triangle commutes. It is enough to check it on generators: let $x
\in VX$ and $s \in G_x$. We have $\Phi(s)=\gamma_x \phi_x(s) {\gamma_x}^{-1}$, $\Psi(s)=\gamma'_x \psi_x(s)
{\gamma'_x}^{-1}$ and on the other hand
\begin{equation}
\begin{split}
\Theta \circ \Phi(s)=&\Theta(\gamma_x) \Theta(\phi_x(s))
\Theta(\gamma_x)^{-1}\\
=&\Theta(\gamma_x)\rho_{\phi(x)}\theta_{\phi(x)}(\phi_x(s))\rho_{\phi(x)}^{-1}
\Theta(\gamma_x)^{-1}\\
=&\Theta(\gamma_x)\rho_{\phi(x)}\zeta_x\psi_x(s)\zeta_x^{-1}\rho_{\phi(x)}^{-1}
\Theta(\gamma_x)^{-1}\\
\end{split}
\end{equation}
(using property (b) of isomorphism of coverings), and this is equal to
$$=\gamma'_x \psi_x(s) {\gamma'_x}^{-1}=\Psi(s)$$
by property (c) and the definition of $\Psi$. Similarly, for $e \in EX$,
\begin{equation}
\begin{split}
\Theta \circ \Phi(e)=&\Theta(\gamma_e) \Theta(\phi(e))
\Theta(\gamma_{\bar{e}})^{-1}\\
=&\Theta(\gamma_e)\rho_{\phi(e)}\theta(\phi(e))\rho_{\phi(\bar{e})}^{-1}
\Theta(\gamma_{\bar{e}})^{-1}\\
=&\Theta(\gamma_e)\rho_{\phi(e)}\psi(e)\rho_{\phi(\bar{e})}^{-1}
\Theta(\gamma_e)^{-1}= \gamma'_e \psi(e) {\gamma'_{\bar{e}}}^{-1}.\\
\end{split}
\end{equation}
The last equality comes from the fact that since $\zeta_e \in H'_{\psi(e)}$, by definition of the fundamental
group,$$\psi(e)=\zeta_e \psi(e) \zeta_{\bar{e}}^{-1}.$$
 Thus we have a commuting triangle of morphisms of fundamental groups
$$\xymatrix{
\pi_1(X,G_{\bullet},x_0) \ar[r]^-{\Phi} \ar[dr]^-{\Psi}& \pi_1(Y,H,\phi(x_0)) \ar[d]^-{\Theta}\\
& \pi_1(Y',H',\psi(x_0)) }$$ (where $\Theta$ is an isomorphism). In a similar fashion, we have a triangle of
\textit{isomorphisms} of trees, which is equivariant with respect to the above triangle of groups:

$$\xymatrix{
(\widetilde{X,G_\bullet,x_0}) \ar[r]^-{\tilde{\phi}} \ar[dr]^-{\tilde{\psi}}&
(\widetilde{Y,H_\bullet,\phi(x_0)})
\ar[d]^-{\tilde{\theta}}\\
& (\widetilde{(Y',H',\psi(x_0)}) }.$$ We claim that this triangle is also commutative. Indeed, by definition, if
$g \in \pi[x_0,x]/G_x \subset (\widetilde{X,G_{\bullet},x_0})$ then
\begin{equation*}
\begin{split}
\tilde{\theta}(\tilde{\phi}(g))&=\tilde{\theta}(\Phi(g)\gamma_x)
=\Theta(\Phi(g))\Theta(\gamma_x)\rho_{\phi(x)}\\
&=\Psi(g)\gamma'_x \zeta_x^{-1} =\Psi(g)\gamma'_x=\tilde{\psi}(g)
\end{split}
\end{equation*}
where we used relation (c) together with the fact that $\zeta_x\in H'_{\psi(x)}$ (observe that $\Psi(g)\gamma'_x
\in \pi[\psi(x_0),\psi(x)]/H'_{\psi(x)}$).

\paragraph{}In a similar way as we defined $\Gamma_{\phi}$, define
$\Gamma_{Y'} \subset Aut((\widetilde{Y',H'_{\bullet},\psi(x_0)}))$, an embedding $i_{\psi}: \Gamma_{Y'} \to
Aut((\widetilde{X,G_{\bullet},x_0}))$ and put $\Gamma_{\psi}=Im(i_\psi) \subset
Aut((\widetilde{X,G_{\bullet},x_0}))$. We claim that $\Gamma_{\phi}=\Gamma_{\psi}$. Indeed, if $u \in \Gamma_Y$
then $\Theta(u) \in \Gamma_{Y'}$ and for $v \in (\widetilde{X,G_{\bullet},x_0})$ we have
\begin{equation*}
\begin{split}
i_\phi(u)\cdot v &= \tilde{\phi}^{-1}(u \cdot \tilde{\phi}(v))
=\tilde{\phi}^{-1}(\tilde{\theta}^{-1}(\Theta(u)
\cdot
\tilde{\theta}(\tilde{\phi}(v)))\\
&=\tilde{\psi}^{-1}(\Theta(u) \cdot \tilde{\psi}(v))
 =i_{\psi}(\Theta(u)) \cdot v.
\end{split}
\end{equation*}
We deduce that $\Gamma_{\phi} \subset \Gamma_{\psi}$. Replacing $\theta_{\bullet}$ by its inverse and exchanging
the roles of $\psi_{\bullet}$ and $\phi_{\bullet}$ we obtain the reverse inclusion $\Gamma_{\psi} \subset
\Gamma_{\phi}$. Thus $\Gamma_{\psi}=\Gamma_{\phi}$ as desired.\qed

\medskip
\textit{Proof of Proposition 1.5.} It remains to show that the map $\phi_\bullet \mapsto \Gamma_\phi$ is the
inverse map of $m$. To see this, let $\Gamma' \supset \Gamma$ be an overlattice of $\Gamma$. The quotient graph
of groups $\Gamma \backslash \backslash T=(X,G_\bullet)$ is formed relative to some datum $(R,S,g_x)$; let us
similarly choose datum $(R',S',g'_x)$ inducing a quotient graph of groups $(Y,H_\bullet)=\Gamma' \backslash
\backslash T$. Recall that by \cite{S}, \S 5.4, there are, for any $x_0 \in VX$ and $y_0 \in VY$, canonical
isomorphisms $\Gamma \simeq \pi_1(X,G_\bullet,x_0)$, $T \simeq (\widetilde{X,G_\bullet,x_0})$ and $\Gamma'
\simeq \pi_1(Y,H_\bullet,y_0)$, $T \simeq (\widetilde{Y,H_\bullet,y_0})$. Choosing furthermore some elements
$\theta_x$ as in the proof of Lemma \ref{L:1} we get a covering (see \cite{B}, Section 4.2)
$$m^{\Gamma'}:\; (X,G_\bullet) \to (Y,H_\bullet).$$
>From \cite{B}, Proposition 4.2, the following diagrams commute :
$$\xymatrix{
T \ar[d] \ar[r]^{id} & T \ar[d]\\
(\widetilde{X,G_\bullet,x_0}) \ar[r]^{\widetilde{m^{\Gamma'}}} & (\widetilde{Y,H_\bullet,y_0})}$$

$$\xymatrix{
\Gamma \ar[d] \ar[r] & \Gamma' \ar[d]\\
\pi_1(X,G_\bullet,x_0) \ar[r]^{M^{\Gamma'}} & \pi_1(Y,H_\bullet,y_0)}$$ where we denote $M^{\Gamma'}$ the
morphism of path groups induced by the covering $m^{\Gamma'}$.

In particular, the pullback of $\pi_1(Y,H_\bullet,y_0)$ via the composition of isomorphisms $T \simeq
(\widetilde{X,G_\bullet,x_0}) \stackrel{\widetilde{m^{\Gamma'}}}{\to} (\widetilde{Y,H_\bullet,y_0})$ is equal to
$\Gamma'$. This shows that $\phi_\bullet \mapsto \Gamma_\phi$ is a left inverse of $\Gamma' \mapsto
m^{\Gamma'}$. The other direction is proved in a similar way.

 Let $\Gamma'$ be an overlattice of $\Gamma$ of index $n$. Note that $m^\Gamma$ is an $n'$-sheeted covering with $n=n'$,
 as we have
\begin{align*}
n=[\Gamma': \Gamma] &= \frac{ vol (\Gamma\backslash \backslash T)}{ vol (\Gamma' \backslash \backslash T)}=
\frac{\underset{x \in VX}{\sum} \frac{1}{|G_x|}}{\underset{y \in VY}{\sum} \frac{1}{|H_y|}} =\frac{\underset{y
\in VY}{\sum} \underset{x \in \phi^{-1}(y)}{\sum} \frac{1}{|G_x|}}{\underset{y \in VY}{\sum} \frac{1}{|H_y|}}=
\frac{\underset{y \in VY}{\sum} \frac{n'}{|G_y|}} {\underset{y \in VY}{\sum} \frac{1}{|H_y|}} =n'.
\end{align*}

Note that the first equality comes from the fact that $T$ is a left $\Gamma'$-set (and $\Gamma$-set) with finite
stabilizers (see \cite{BL}, page 16). \qed

\paragraph{}
It follows from the above proposition  that to find $u(n)$, it suffices to count the number of isomorphism
classes of coverings of faithful graphs of groups by $(X,G_\bullet)$.

\section{Main results}
\paragraph{2.1}
Let $G$ be a group of order $n$ and let $n=\prod_{i=1} ^t p_i ^{k_i}$ be the prime decomposition of $n$. Let
$\mu=\mu(n)$ be the maximum of $k_i$.
We denote by $d(G)$ the minimal cardinality of a generating set of $G$ and by
$f(n)$ the number of isomorphism classes of groups of order $n$.

In \cite{P}, Pyber showed that the number of isomorphism classes of groups of order $n$ with a given Sylow set,
namely the set of Sylow $p_i$-subgroups defined up to conjugacy, is at most $n^{75\mu+16}$. Together with the
result of Sims (\cite{Si}), namely $f(p^k) \leq p^{\frac{2}{27} k^3 + \frac{1}{2} k^{\frac{8}{3}}}$, we get the
following upper bound for $f(n)$:
\begin{align*}
f(n) &\leq \prod_{i=1}^{t} p_{i}^{\frac{2}{27}k_{i}^3 + \frac{1}{2}k_i^{8/3}}n^{75\mu+16}\\
&\leq n^{\frac{2}{27}\mu^2 +\frac{1}{2}\mu^{5/3}+75\mu+16}
\end{align*}
Let $g(n)=\frac{2}{27}\mu^2(n) + \frac{1}{2}\mu^{5/3}(n)+ 75\mu(n) +16$ so that $f(n) \leq n^{g(n)}$.

On the other hand, Lucchini and Guralnick showed that if every Sylow subgroup of G can be generated by $d$ elements,
then $d(G) \leq d+1$ (\cite{Luc}, \cite{Gu}). Combining with the basic fact that  $d(H) \leq n$ for any group $H$ of
order $p^n$ (\cite{Si}), we deduce that
\begin{equation*}
d(G) \leq \mu +1.
\end{equation*}

 Using these results, we obtain the following upper-bound for $u(n)$.

\begin{theo} Let $\Gamma$ be a cocompact lattice of $Aut(T)$. Then there are some positive constants $C_0$ and $C_1$
depending only on $\Gamma$, such that
\begin{equation*}
\forall n > 1, \;\;\;\;\;\; u_{\Gamma} (n) \leq C_0 n^{C_1 log^2(n)}.
\end{equation*}
\end{theo}
\begin{lem} Any covering $\phi_\bullet=(\phi,\phi_x,\gamma_x) : (X,G_{\bullet}) \to (Y,H_{\bullet})$ is isomorphic to a
covering $\phi'_\bullet=(\phi',\phi'_x,\gamma'_x) : (X,G_\bullet) \to (Y',H'_\bullet)$ where $\gamma'_x \in
 \Pi(Y',H'_\bullet)$ is a product of at most $12 \times \text{diameter}(X)$ generators $h_y \in H'_y$ and $e \in EY'$.
 \end{lem}

\noindent \textit{Proof.} Fix $x_0 \in X$. Associated to $\phi_\bullet$ is a lattice $\Gamma' \subset
Aut((\widetilde{X,G_\bullet, x_0}))$ containing $\pi_1(X,G_\bullet,x_0)$.  From $(X,G_\bullet,x_0)$ we construct
$(R,S,g_e)$ such that the quotient of $(\widetilde{X,G_\bullet,x_0})$ by $\pi_1(X,G_\bullet,x_0)$ is exactly
$(X,G_\bullet)$. Namely, first fix a maximal
 tree $\tau$ in $X$. We may choose
 $R=\{e_1\cdots e_n\;|\;e_1\cdots e_n\;\\\text{is\;a\;path\; from}\; x_0 \; \text{in} \tau\}$, $S= \{e_1\cdots e_n e_{n+1}\;|
 \;e_1\cdots e_n\;
 \text{is\;a\;path\;in\;}\tau\}$ and $g_e=e'_1 \cdots e'_l {e_{n+1}}^{-1} \cdots {e_1}^{-1}$ where $e'_1
 \cdots e'_l$ is a path in $\tau$ from $x_0$ to $t(e)$, and where $e$ is the edge connecting $e_1 \cdots e_n$ to
 $e_1 \cdots e_{n+1}$. In particular, $g_e$ is a product of at most twice the diameter of $X$ number of generators of $\Pi(X,G_\bullet)$. Now we choose
 $R', S'$ subsets of $R,S$ in such a way that the restriction of the projection $(\widetilde{X,G_\bullet,x_0}) \to \Gamma'
 \backslash (\widetilde{X,G_\bullet,x_0})$ on $R'$ is bijective for vertices (resp. the restriction of $S'$ is bijective
 on edges). We also choose $g'_e$ in a similar fashion as above, hence $g'_e$ is also a product of at most twice
the diameter of $X$ number of generators of $\Pi(X,G_\bullet)$. From this data, we construct a graph of groups
$(Y',H'_\bullet)$ as usual, and we have a canonical injection $\Gamma' \subset \Pi(Y',H'_\bullet)$. For $x \in
X$ there exists a unique lift $\tilde{x} \in R$ and a unique $\tilde{x'} \in R'$ in the $\Gamma'-$orbit of
$\tilde{x}$. Choose $\theta_x \in \Gamma' \subset \Pi(X,G_\bullet)$ such that $\theta_x(\tilde{x})=\tilde{x'}$
and such that $\theta_x$ is a product of at most $l(\tilde{x},\tilde{x_0}) + l(\tilde{x_0},\tilde{x'})$
generators of $\Pi(X,G_\bullet)$: here $l(a,b)$ is the distance in the tree $R$ between $a$
 and $b$. This is possible since we may first choose a path in the path group $\Pi(X,G_\bullet)$
 from $\tilde{x}$ to $\tilde{x_0}$ of length $\leq l(\tilde{x},\tilde{x_0})$ and then a path
 from $\tilde{x_0}$ to $\tilde{x'}$ of length $\leq l(\tilde{x_0},\tilde{x'})$. Observe that since we chose $R'
 \subset R$, we have $l(\tilde{y},\tilde{w}) \leq \mathrm{diameter}(X)$ for any vertices $w,y \in VX$.
 We do the same thing for edges in $S$, to define $\theta_e \in
 \Pi_1(X,G_\bullet)$ such that $\theta_e(\tilde{e})=\tilde{e'}$
and $\theta_e$ is a product of at most $2 \times \text{diameter}(X)$ generators of $\Pi(X,G_{\bullet})$. Then we
can construct from $\theta_x$ and $\theta_e$'s a covering $\phi'_\bullet: (X,G_\bullet) \to (Y',H'_\bullet)$,
with $\gamma'_x=\theta_x^{-1}$ and $\gamma'_e=g_{e}\theta_e^{-1}g_{e'}^{-1}$, which are both products of at most
$6 \times \text{diameter}(X)$ generators of $\Pi(X,G_\bullet)$. Observe that a word of length $l$ in generators
of $\Pi(X,G_\bullet)$ belonging to $\Gamma'$ is also expressible as a word of length $l$ in generators of
$\Pi(Y',H'_\bullet)$.
 Finally, by the proposition on bijection of isomorphism
classes of coverings and overlattices, $\phi_\bullet: (X,G_\bullet) \to (Y,H_\bullet)$ is isomorphic to
$\phi'_\bullet$. \qed

\paragraph{}
\noindent \textit{Proof of Theorem 2.1.} Let us fix a quotient graph of groups $(X,G_\bullet)$ of $\Gamma$ as in
section 1.1. There exist only finitely many coverings of edge-indexed graphs by the edge-indexed graphs
underlying $(X, G_\bullet)$, thus it is enough to show the assertion for the number of overlattices with a fixed
edge-indexed graph. Thus we want to count n-sheeted covering graphs of groups $\phi_\bullet:(X,G_\bullet) \to
(Y, H_\bullet)$ such that $Y$ is a fixed subgraph (with fixed indices) of $X$  and $\phi : X \rightarrow Y$ the
natural projection, and that the edge group $H_e$ is a subgroup of $H_{o(e)}$.

 Let $c_x = |G_x|$ for any $x$ in $VX \cup EX$ and let $c_y = (\underset{x \in \phi^{-1}(y)}{\sum} c_x^{-1})^{-1}$.
 By the definition of $n$-sheeted covering, the cardinality $|H_y|=n c_y$, for any $y$ in $VY \cup EY$.

Now we claim that for any group $H$ of order $n$, there are at most $(m!)^{\mu(n)+1}$   subgroups of index $m$.
For to any transitive $H$-action on the set $\{1, \cdots, m \}$, we can associate a subgroup of $H$ with index
$m$, namely the stabilizer of $1$. This map $\{\rho: H \to S_m\} \longrightarrow \{H' \subset H | [H:H']=m \}$
is surjective since for any subgroup $H'$ of $H$ with index $m$, the action of $H$ on the cosets $H/H'$ gives
(among many) an action on $\{1, \cdots, m \}$, where we let $1$ stand for the trivial coset $H'$. Again by the
theorem of Luccini and Guralnick, there are at most $(m!)^{\mu(n)+1}$ transitive $H$-action on the set $\{1,
\cdots, m \}$, as claimed.

 There are at most $\prod_{y \in VY}(c_y n)^{g(c_y n)}$ isomorphism classes of $H_y$'s. By the above claim, the
 number of subgroups $\alpha_f(H_f)$ of $H_y$ is at most $((c_y/c_f)!)^{\mu(c_y n)+1}$. There are at most
 $\prod_{f \in EY} (c_f n)^{\mu(c_f n)+1}$ isomorphisms $\varphi: \alpha_f H_f \to \alpha_{\bar{f}}H_f$ and at
 most $\prod_{x \in VX} (c_{\phi(x)} n)^{\mu(c_x)+1}$ injections $\phi_x :G_x \to H_{\phi(x)}$.
By Lemma 2.2, there are at most $(max_y |H_y|)^{12K}$ choices for each $\gamma_x$ or $\gamma_e$, where
$K=diameter\;of\;X$. Hence
$$\#\{(\gamma_x,\gamma_e)\} \leq \prod_{x \in VX} max_{y \in VY} (c_yn)^{12K} \times \prod_{e \in EX} max_{y \in VY}
(c_yn)^{12K}$$ which is bounded by $(Mn)^{(12K)(|VX|+ |EX|)}$.

 Note that by the condition of injectivity and the commutativity of the diagram,

$$\xymatrix{
G_e \ar[d]^{\phi_e} \ar[r]^{\alpha_e} & G_x \ar[d]^{\phi_x}\\
H_{\phi(e)} \ar[r]^-{ad(\gamma_x^{-1} \gamma_e) \circ \alpha_{\phi(e)}} & H_{\phi(x)}}$$

the group morphism $\phi_e:G_e \to H_\phi(e)$ is completely determined by the morphism $\phi_x:G_x \to H_x$.

Let $M=\underset{y \in VY \cup EY} \max c_y$, $\mu=\mu(Mn)$. Let $c_0=|VY|$, $c_1=|EY|$, $c_2= \underset{\{f \in
EY\}} {max} \lbrace\big(\frac{c_{o(f)}}{c_f} \big)!\rbrace$, let $c_3= \underset{x \in VX}{\sum} \mu(c_x)+1$.
Combining all the estimates above, we get the following upper bound for $u(n)$,
\begin{align*}
u_\Gamma(n)&\leq \prod_{y \in VY}(c_y n)^{g(c_y n)}
                 \prod_{x \in VX}(c_{\phi(x)} n)^{\mu(c_x)+1}
                 \prod_{f \in EY} (c_f n)^{\mu(c_f n)+1}
                 \prod_{f \in EY}((c_{o(e)}/c_e)!)^{\mu_{c_{o(e)} n}+1}
         \cdot ((Mn)^{12K(|VX|+|EX|)}) \\
           &\leq \prod_{y \in VY}(M n)^{g(M n)}
                 \prod_{x \in VX}(M n)^{\mu(c_x)+1}
                 \prod_{f \in EY}(M n)^{\mu(M n)+1}
                 \prod_{f \in EY}(c_2)^{\mu_{M n}+1}
         \cdot ((Mn)^{12K(|VX|+|EX|)}) \\
           &\leq (Mn)^{c_0g(Mn)+ c_3 + c_1(\mu(Mn)+1)}
                 (c_2)^{c_1(\mu(Mn)+1)} (Mn^{12K(|VX|+|EX|)}) \\
           &\leq (Mn)^{\frac{2}{27}c_0\mu^2 + \frac{c_0}{2}\mu^{5/3}+
                       (75c_0+2c_1) \mu +(16c_0+2c_1+c_3)+c_2^{c_1(\mu+1)}}
                       ((Mn)^{12K(|VX|+|EX|)})\\
           &\leq (C_0 n)^{C_1 \mu^2} \leq (C_0 n)^{C_1' (\log n)^2} \\
\end{align*}
where $C_0=max\{M, c_2\}\cdot (Mn^{12K(|VX|+|EX|)})$, $C_1= c_0\big(\frac{2}{27}+\frac{1}{2} +75+16+2
\big)+3c_1+c_3+12K(|VX|+|EX|)$ and $C_1'=\frac {C_1}{(\log 2)^2}$. \qed

\paragraph{2.2.}
Let $p$ be a prime number.
>From now on, we assume that $T$ is a $2p$-regular tree and that $\Gamma$ is a cocompact lattice in $Aut(T)$ with a
quotient graph of groups given by
$$
\centerline{\epsfbox{essai.eps}}
$$


The aim of this section is to give, in this situation, a smaller upper bound on $u_\Gamma(n)$ than the previous one,
as well as a lower bound.

\begin{theo}
Let $n=p_0 ^{k_0} p_1 ^{k_1} \cdots p_t ^{k_t}$ be the prime decomposition of $n$ with $p_0=p$. Then there exist
positive constants $c_0, c_1$ such that $\underset{k_0 \to \infty}{\limsup} \frac{u_\Gamma(n)}{n^{c_1 \log n}} \leq c_0$.
For $n=p_0^{k_0} (k_0 \geq 3)$, we also have $u_\Gamma(n)\geq n^{\frac{1}{2}(k-3)}$.
\end{theo}

 In the following lemma, we denote by $[g,h]$ the commutator $ghg^{-1}h^{-1}$ in $G$.

\begin{lem} Let $A=(a_{s,t})_{1 \leq s,t \leq k-1}$ be a lower triangular matrix with coefficients in ${0,\cdots, p-1}$
and $G=G(A)$ be a group defined by the generators $\bar{g}_0, \bar{g}_1, \cdots, \bar{g}_k$ and the following relators
\begin{align*}
 \bar{g}_i ^p &=1,  \;\;\;\;\;\;\;\;\;\;\;\;\;\;\;\;\;\;\;\;  i=0,1, \cdots, k\\
[\bar{g}_i, \bar{g}_{i+1}] &=1, \;\;\;\;\;\;\;\;\;\;\;\;\;\;\;\;\;\;\;\; i=0, 1, \cdots, k-1  \\
[\bar{g}_i, \bar{g}_{i+2}] &=\bar{g}_{i+1}^{a_{1,1}}, \;\;\;\;\;\;\;\;\;\;\;\;\;\;\;  i=0, 1, \cdots, k-2
\;\;\;\;\;\;\; \tag{**} \\
[\bar{g}_i, \bar{g}_{i+3}] &=\bar{g}_{i+1}^{a_{2,1}}\bar{g}_{i+2}^{a_{2,2}}, \;\;\;\;\;\;\;\; i=0, 1, \cdots, k-3\\
       &\vdots\\
[\bar{g}_i, \bar{g}_{i+k}] &=\bar{g}_{i+1}^{a_{k-1,1}}\bar{g}_{i+2}^{a_{k-1,2}}\cdots \bar{g}_{i+k-1}^{a_{k-1,k-1}},
\;\;\;\;\; i=0
\end{align*}
 Then any element of $G$ can be written as $\bar{g}_0^{i_0}\cdots \bar{g}_k^{i_k}$ where $0 \leq i_j < p$.

\end{lem}
\noindent \textit{Proof.} We proceed by induction on $k \geq 1$. It is clear for the case $k=1$ since the
generators $\bar{g}_0$ and $\bar{g}_1$ commute. Now suppose that the assertion is true for all $k \leq m-1$. For
$k=m$, consider the subgroup $G_1$ generated by $\bar{g}_0, \bar{g}_1, \cdots, \bar{g}_{m-1}$. It is a quotient
of $G(A')$ with $A'=(a_{s,t})_{1 \leq s, t \leq k-2}$, by induction hypothesis, any element of $G_1$ can be
written as $\bar{g}_0^{i_0}\cdots \bar{g}_{m-1}^{i_{m-1}}$ where $0 \leq i_j < p$. Now we only need to consider
the elements of $G-G_1$. By an easy induction, it suffices to consider the elements $\bar{g}_m
\bar{g}_i=[\bar{g}_i,\bar{g}_m]^{-1}\bar{g}_i\bar{g}_m$ for $i=0, 1, \cdots, m-1$. Since $[\bar{g}_i, \bar{g}_m]
\in [G, G] \subset G_1$, the element $[\bar{g}_i, \bar{g}_m]^{-1}\bar{g}_i$ is an element in $G_1$, thus can be
expressed as $\bar{g}_0^{i_0}\bar{g}_1^{i_1}\cdots \bar{g}_{m-1}^{i_{m-1}}$ for some $i_j$ in $\{0, \cdots,
p-1\}$. Therefore we get $\bar{g}_m \bar{g}_i =\bar{g}_0^{i_0} \cdots \bar{g}_{m-1}^{i_{m-1}}\bar{g}_m$.
\;\;\;\;\;\qed

\begin{lem} Let $G$ be a group of order $p^{k+1} (k \geq 1)$, $G_{1}$ and $G_{2}$ be two isomorphic subgroups of index
$p$ in $G$ and $\varphi$ be an isomorphim from $G_{1}$ to $G_{2}$. Suppose that $G_1$ contains no subgroup $N$
which is normal in G and $\varphi$-invariant. Then
\begin{enumerate}
\item[(a)]there exist elements $g_{i}$ in $G$ for $i=0, \cdots , k$, such that  $\varphi(g_{i})=g_{i+1}$ for $i=0,
\cdots , k-1$, $G=\langle g_{0},\cdots , g_{k}\rangle$ and $G_{1}=\langle g_{0},\cdots , g_{k-1}\rangle$,
\item[(b)]There exists a lower triangular matrise $A$ with coefficients in $0, \cdots, p-1$ such that the map
$\psi : G(A) \to G$ defined by $\bar{g}_i \mapsto g_i$ is well-defined and is an isomorphism.

\end{enumerate}
\end{lem}
\noindent \textit{Proof.} We proceed by  induction on $k \geq 1$, using the fact that any maximal proper
subgroup of a $p$-group is normal, see for instance \cite{Su}.

We first consider the case $k=1$, that is when $G$ has order $p^2$. Since $G_1$ and $G_2$ are maximal, they are
normal. Thus they are not equal by the normality assumption and $G=\langle G_1, G_2\rangle $. Let $g_0$ be an
element in $G_1 - G_2$ and set $g_1 =\varphi(g_0)$. Then clearly $G_1=\langle g_0\rangle , G_2=\langle
g_1\rangle $ and $G=\langle g_0, g_1\rangle $. Moreover, since $|G|=p^2$, $G$ is abelian \cite{Su}. Thus $[g_0,
g_1]=1$ and $G=\{g_0^{i_0} g_1^{i_1}: 0 \leq i_0, i_1 \leq p-1\}$, which shows that $\psi$ in (b) is
well-defined and surjective. Since $G(A)$ has cardinality at most $p^2$ by Lemma 2.3, and $G$ has cardinality
$p^2$, the map $\psi$ is an isomorphism.

Now suppose that the assertion is true for all $k < m$. For $k=m$, consider $G, G_1, G_2$ and $\varphi$ as in the
statement of the lemma.
As above, $G_1$ and $G_2$ are normal, distinct and $G=\langle G_1, G_2\rangle$. Since $[G_2 : G_2 \cap G_1]\leq
[\langle G_1, G_2\rangle :G_1]=p$, we have $ [G_2:G_1 \cap G_2]=p$ and similarly $[G_1:G_1 \cap G_2]=p$. Therefore
$G_1 \cap G_2$ is maximal, thus normal in $G_1$ and $G_2$. Since $G$ is generated by $G_1$ and $G_2$, the subgroup
$G_1 \cap G_2$ is normal in $G$. By the assumption, $\varphi(G_1 \cap G_2) \neq G_1 \cap G_2$.\\
\\
\noindent
\textbf{Claim.} If a subgroup $N$ of $G_1 \cap G_2$ is normal in $G_2$ and $\varphi$-invariant, then $N$ is normal in
$G$.\\
$$
\xymatrix @-1.2pc{
&G \ar@{-}[dl] \ar@{-}[dr]\\
G_1 \ar@{-}[dr] \ar[rr]^{\varphi} & & G_2 \ar@{-}[dl] \ar@{-}[dr] \\
&{G_1 \cap G_2}  \ar@{-}[dr] \ar[rr]^{\varphi} && \varphi(G_1 \cap G_2) \ar@{-}[dl]\\
&&N
}
$$
\noindent \textit{Proof.} Consider $gNg^{-1}$, for any $g \in G_1$. As $\varphi(gNg^{-1})
=\varphi(g)\varphi(N)\varphi(g)^{-1} = \varphi(g)N \varphi(g)^{-1}=N$ (since $\varphi(g) \in G_2$) and $\varphi$
is an isomorphism, we deduce that $gNg^{-1}=\varphi^{-1}(N)=N$. Therefore $G_1 \subseteq N_{G} (N)$ and
similarly  $G_2 \subseteq N_{G} (N)$. Thus $G$, as a group generated by $G_1$ and $G_2$, is also contained in $
N_{G} (N)$. Hence $N$ is normal in $G$. \;\;\;\;\qed
\\
\noindent By the above claim, we can use the induction hypothesis on $G'=G_2$, $G_1'=G_1 \cap G_2$,
$G_2'=\varphi(G_1 \cap G_2)$ and $\varphi'=\varphi|_{G_1 \cap G_2}$. It follows that there is an element $g_1$
in $G_2$ such that $G_2 =\langle g_1, \cdots, g_m\rangle$, $G_1 \cap G_2 =\langle g_1, \cdots, g_{m-1}\rangle $
and $\varphi(g_i)=g_{i+1}$ for $i=1, \cdots , m-1$.

Let $g_0 =\varphi^{-1} (g_1)$. If $g_0 \in G_1 \cap G_2$, then $g_1 \in \varphi(G_1 \cap G_2)$, which contradicts
$G_1 \cap G_2 \neq \varphi(G_1 \cap G_2)$. Thus $g_0$ is an element of $G_1 - G_2$. Since $G_2$ is maximal in $G$,
the group $G$ is generated by $G_2$ and $g_0$, i.e., $G=\langle g_0, g_1, \cdots, g_m\rangle $. It is clear that
$G_1=\langle g_0, \cdots, g_{m-1}\rangle$ as $G_1 \cap G_2$ is maximal in $G_1$ and $g_0 \in G_1-(G_1 \cap G_2)$.

To prove the assertion (b), note that $[g_0, g_i]=\varphi^{-1}([g_1, g_{i+1}])=\varphi^{-1}(g_2^{a_{i,1}}\cdots
g_i^{a_{i,i}})=g_1^{a_{i,1}}\cdots g_{i-1}^{a_{i,i}}$ for $1 \leq i \leq m-1$ and ${g_i}^p=1$, for all $i \geq
1$ by induction hypothesis. Thus we only need to consider $[g_0, g_m]$ and ${g_0}^p$. The element $g_0$ clearly
has order $p$ since $\varphi$ is an isomorphism and $g_1=\varphi(g_0)$ has order $p$. It is easy to see that if
two subgroups $H$ and $K$ are normal subgroups of a group $G$, then so is the commutator subgroup $[H,K]$ and we
have $[H, K] \subset H \cap K$. Since $g_0 \in G_1$ and $g_m \in G_2$, it follows that $[g_0, g_m] \in G_1 \cap
G_2 =\langle g_{1}, g_{2}, \cdots, g_{m-1}\rangle,$ which proves that $\psi$ is a well-defined homomorphism. By
the previous paragraph, it is surjective. Since $G$ and $G(A)$ are of cardinality $p^{k+1}$ and at most
$p^{k+1}$ respectively, the map $\psi$ is an isomorphism. \qed
\paragraph{}
\noindent
\textit{Proof of Theorem 2.2.} By Proposition 1.2, the number $u_n$ is the number of isomorphism classes of
$n$-sheeted  coverings of faithful graphs of groups $\phi_\bullet:\Gamma \backslash \backslash T \to (X, G_\bullet)$.
As already seen, we may assume that $X=\Gamma \backslash T$. The following commutative diagram summarizes the data
defining $\phi_\bullet$:
$$
\xymatrix{
  &1 \ar[d] \ar[r] &\mathbb{Z}/p\mathbb{Z} \ar[d] \\
  &G_e \ar@<2pt>[r]^{\alpha_{e}} \ar@<-2pt>[r]_{\alpha_{\overline{e}}}  &G_x}
$$

 Let's first consider the case when $n=p^k$. Let $G=G_x$, $G_1 =\alpha_e (G_e)$ and  $G_2 = \alpha_{\overline{e}}
(G_e)$. By the condition of faithfulness, $G_1$ and $G_2$ are distinct as they are normal subgroups of $G$.
Hence if we let $\varphi= \alpha_{\bar{e}} \circ \alpha_e^{-1} :G_1 \to G_2$, then $\varphi$ is an isomorphism
and there is no subgroup of $G_1$ which is normal in $G$ and $\varphi$-invariant. Thus we can use Lemma 2.4 to
find an element $g_0$ in $G_x$ such that $G_x=\langle g_0, g_1, \cdots g_k \rangle$ where
$\varphi(g_j)=g_{j+1}$.
 Moreover, the group $G$ is isomorphic to $G(A)$, which is determined by $A$. (note that $A$ also determines $G_e$
and the maps $\alpha_e$ and $\alpha_{\bar{e}}$.) Thus we have at most $p^{\sum_{j=0}^{k-1}j}$ choices for $G_x,
G_e, \alpha_e$ and $\alpha_{\bar{e}}$, which is exactly the number of choices of $(a_{st})_{1 \leq t \leq s \leq
k-1}$. Once we have fixed $G_x$, $G_e$, $\alpha_e$ and $\alpha_{\bar{e}}$, an injection $i$ from
$\mathbb{Z}/p\mathbb{Z}$ into $G_x$ is determined by the image of a generator in the domain, which implies that
we have at most $|G_x|=p^{k+1}$ choices for $i$. Therefore we have an upper bound $u_\Gamma(n) \leq
p^{\sum_{j=1}^{k-1}j} p^{k+1}=p^{\frac{(k-1)k}{2}+k+1}=p^{\frac{k^2+k+2}{2}}$.

Now let us construct non-isomorphic classes of  faithful covering graph of groups to deduce a lower bound of
$u(n)$. Fix a lower triangular matrix $A=(a_{st})$ $1 \leq t \leq s \leq k-1$ with coefficients in ${0, \cdots,
p-1}$ such that furthermore $a_{k-1,j}=0$ for $j=0, \cdots , k-1$. Let $G_x=G(A)$ be the group defined in Lemma
2.3, and define the  covering graph of groups $\phi_\bullet=\phi_\bullet(A): \Gamma \backslash \backslash T \to
G_\bullet$ as follows. Let $G_e$ be the subgroup of $G_x$ generated by $\bar{g}_0, \cdots, \bar{g}_{k-1}$. Let
the injection $\alpha_e$ be the inclusion map and the other inclusion $\alpha_{\bar{e}}$ be defined by
$\alpha_{\bar{e}}(\bar{g}_i)=\bar{g}_{i+1} \in G_x$, which is indeed a monomorphism by the definition of $G(A)$.
Therefore, the group morphism $\varphi=\varphi_{A}$ defined by $\varphi(\bar{g}_i)=\bar{g}_{i+1}$ is an
isomorphism from $\alpha_e(G_e)$ onto $\alpha_{\bar{e}}(G_e)$. The data $G_\bullet$ thus defines a faithful
graph of groups, as there is no $\varphi$-invariant subgroup of $G_e$. Indeed, for any nontrivial element
$h=\bar{g}_{j_0} ^{i_{j_0}}\cdots \bar{g}_{j_t} ^{i_{j_t}}$ in $G_e$ (with nonzero $i_{j_t}$),
$\varphi^{k-j_t}(h)= \bar{g}_{j_0+k-j_t}^{i_{j_0}}\cdots \bar{g}_k ^{i_{j_t}} \notin G_e$.
 Let $v$ be a generator of $\mathbb{Z}/p\mathbb{Z}$ and set $\phi_x(v)=\bar{g}_0 \bar{g}_k$. This defines a group
monomorphism $\phi_x:\mathbb{Z}/p\mathbb{Z} \to G_x$, as $\phi_x(v)^p=(\bar{g}_0\bar{g}_k)^{p}=\bar{g}_{0}^p
\bar{g}_{k}^p=1$. The map $\phi_x$ is clearly injective since the order of $g_0 g_k$ is $p$.  Thus we have
constructed a  covering of graphs of groups. Now suppose that the  coverings of graphs of groups
$\phi_\bullet(A)$ and $\phi_\bullet(A'): \Gamma \backslash \backslash T \to (X, G'_\bullet)$ are isomorphic. Let
us denote by $\psi:(X, G_\bullet) \to (X, G'_\bullet)$ an isomorphism between them.  Then there exists a
commutative diagram as follows:
\begin{center}
\begin{picture}(0,0)%
\includegraphics{comm.pstex}%
\end{picture}%
\setlength{\unitlength}{1184sp}%
\begingroup\makeatletter\ifx\SetFigFont\undefined%
\gdef\SetFigFont#1#2#3#4#5{%
  \reset@font\fontsize{#1}{#2pt}%
  \fontfamily{#3}\fontseries{#4}\fontshape{#5}%
  \selectfont}%
\fi\endgroup%
\begin{picture}(4962,4038)(1951,-5194)
\put(5326,-4336){\makebox(0,0)[lb]{\smash{\SetFigFont{7}{8.4}{\rmdefault}{\mddefault}{\updefault}{ $\alpha'_e$}%
}}}
\put(2326,-4411){\makebox(0,0)[lb]{\smash{\SetFigFont{7}{8.4}{\rmdefault}{\mddefault}{\updefault}{ $\alpha_e$}%
}}}
\put(3151,-3511){\makebox(0,0)[lb]{\smash{\SetFigFont{7}{8.4}{\rmdefault}{\mddefault}{\updefault}{ $\psi_e$}%
}}}
\put(4801,-5086){\makebox(0,0)[lb]{\smash{\SetFigFont{7}{8.4}{\rmdefault}{\mddefault}{\updefault}{ $\psi_x$ }%
}}}
\put(6676,-4861){\makebox(0,0)[lb]{\smash{\SetFigFont{7}{8.4}{\rmdefault}{\mddefault}{\updefault}{ $G'_x$}%
}}}
\put(5101,-3811){\makebox(0,0)[lb]{\smash{\SetFigFont{7}{8.4}{\rmdefault}{\mddefault}{\updefault}{ $G'_e$}%
}}}
\put(4201,-2311){\makebox(0,0)[lb]{\smash{\SetFigFont{7}{8.4}{\rmdefault}{\mddefault}{\updefault}{  }%
}}}
\put(4876,-2461){\makebox(0,0)[lb]{\smash{\SetFigFont{7}{8.4}{\rmdefault}{\mddefault}{\updefault}{ $\mathbb{Z}/p\mathbb{Z}$}%
}}}
\put(3751,-1411){\makebox(0,0)[lb]{\smash{\SetFigFont{7}{8.4}{\rmdefault}{\mddefault}{\updefault}{ 1}%
}}}
\put(3376,-4861){\makebox(0,0)[lb]{\smash{\SetFigFont{7}{8.4}{\rmdefault}{\mddefault}{\updefault}{ $G_x$}%
}}}
\put(1951,-3736){\makebox(0,0)[lb]{\smash{\SetFigFont{7}{8.4}{\rmdefault}{\mddefault}{\updefault}{ $G_e$}%
}}}
\end{picture}

\end{center}

 Since $\psi_e$ (respectively $\psi_{\bar{e}}$) is a group isomorphism from $G_1=\alpha_e(G_e)$ to
$G'_1=\alpha_e(G'_e)$ (respectively from $G_2=\alpha_{\bar{e}}(G_e)$ to $G'_2=\alpha_{\bar{e}}(G'_e)$), it
follows that $\psi_x$ maps the following hierarchy to the corresponding one in $G'$.
$$
\xymatrix @-1.2pc{
&G \ar@{-}[dl] \ar@{-}[dr]\\
G_1 \ar@{-}[dr] & & G_2 \ar@{-}[dl] \ar@{-}[dr] \\
&{G_1 \cap G_2}  \ar@{-}[dr]  && \varphi(G_1 \cap G_2) \ar@{-}[dl]\\
&&\dots
}
$$
In particular, $\psi_x$ preserves the smallest group in the hierarchy, i.e. $\psi_x(\langle g_k \rangle)=
\langle g'_k \rangle$. Let $b$ be an element in ${0, \cdots, p-1}$ such that $\psi_x(g_k)={g'_{k}}^b$. The map
$\psi_x$ is completely determined since $\psi_x(g_i)={g'_i}^b$ for $i=0, \cdots, k$. In other words, if $A$ and
$bA'$ are not equivalent modulo $p$ for any integer $b= 0, \cdots, p-1$, then $\phi_\bullet(A)$ and
$\phi_\bullet(A')$ are non-isomorphic coverings.  This implies that there are at least
$\frac{p^{\sum_{i=1}^{k-2}i}}{p}$ non-isomorphic  covering graphs of groups given by the above examples
$\phi_\bullet(A)$ with ``non-homothetic'' $A=(a_{i,j})$'s. Therefore we have a lower bound $u_n \geq
\frac{p^{\sum_{i=1}^{k-2}i}}{p} =p^{\frac{k^2-3k}{2}}$.

\paragraph{}Now let's consider the general case. Recall that $|G_e|=\prod_{i=0}^{t} p_i^{k_i}$ and
$|G_x|=\prod_{i=0}^{t} p_i^{k_i} p = p_0^{k_0+1}\prod_{i=1}^{t} p_i^{k_i}$, thus the order of the Sylow
$p_i$-subgroup of $G_e$ and that of $G_x$ are the same for all $i \neq 0$. Let $G_e ^{(p_i)}$ be \textit{a}
Sylow $p_i$-subgroup of $G_e$. For $i \neq 1$, let $G_x ^{(p_i)}= \alpha_e(G_e ^{(p_i)})$. Choose one $p$-Sylow
subgroup $G_x ^{(p)}$ of $G_x$ containing $\alpha_e (G_e ^{(p)})$.

We are now going to show that the faithfulness condition is inherited to the Sylow $p$-subgroups $G_e^{(p)}$,
$G_x^{(p)}$ of $G_e$ and $G_x$, from which we can use the upper bound given in the first part of the proof.
Conjugating $\alpha_{\overline{e}}$ by an element of $G_x$, if necessary, we may assume that
$\alpha_{\overline{e}}(G_e^{(p)}) \subset G_x^{(p)}$, thus we have the following diagram:
$$
\xymatrix{
  &G_e \ar@<2pt>[r]^{\alpha_{e}} \ar@<-2pt>[r]_{\alpha_{\overline{e}}}  &G_x \\
  &G_e^{(p)} \ar@{^{(}->}[u] \ar@<2pt>[r]^{\alpha_{e}} \ar@<-2pt>[r]_{\alpha_{\overline{e}}}  &G_x^{(p)} \ar@{^{(}->}[u]
}
$$
Suppose that $N \vartriangleleft G_e^{(p)}$ and $\mathcal{N}=\alpha_e(N)=\alpha_{\bar{e}}(N) \vartriangleleft
G_x^{(p)}$. Let $\overline{N}=\langle gNg^{-1} : g\in G_e \rangle$ (respectively $\overline{\mathcal{N}}=\langle
g\mathcal{N}g^{-1} : g\in G_x \rangle$) be the smallest normal subgroup of $G_e$ (respectively $G_x$) containing
$N$ (respectively $\mathcal{N}$).

Note that $\alpha_i$ $(i=e, \overline{e})$ induces a bijection between left cosets
\begin{align*}
G_e /G_e^{(p)} &\longrightarrow G_x /G_x^{(p)}\\
gG_e^{(p)} &\longmapsto \alpha_i(g)G_x^{(p)}.
\end{align*}
For since $\alpha_i$ is injective, if $gG_e^{(p)}$ is mapped to $G_x^{(p)}$, then $g$ is in $\alpha_i(G_e) \cap
G_x^{(p)}$, which is a $p$-group in $\alpha_i(G_e)$ containing $\alpha_i(G_e^{(p)})$. Since
$\alpha_i(G_e^{(p)})$ is a Sylow $p$-subgroup, $\alpha_i(G_e^{(p)}) \cap G_x^{(p)}$ is equal to $\alpha_i(G_e)$.
Thus $\alpha_i(g^{-1}h)$ is contained in $G_x^{(p)}$ if and only if $g^{-1}h \in G_e ^{(p)}$ and the map is
injective. It is surjective since the source and the target have the same cardinality. Thus any element $g$ in
$G_x$ can be written as $\alpha_i(g') h_i$ for some $g'\in G_e$, $h_i\in G_x^{(p)}$ and we have
$g\mathcal{N}g^{-1}=\alpha_i(g')h_i\mathcal{N}h_i^{-1} \alpha_i ({g'}^{-1})=\alpha_i(g')\mathcal{N}\alpha_i
({g'}^{-1})$ ($i=1,2$). Therefore $\alpha_1(\overline{N})=\alpha_2(\overline{N})=\overline{\mathcal{N}}$ and it
is normal in $G_x$. As a consequence, $G_e^{(p)}$ and $G_x^{(p)}$ satisfy the condition of faithfulness, i.e.,
there is no subgroup $N$ of $G_e^{(p)}$ such that $\alpha_e(N)=\alpha_{\bar e}(N)$ is normal in $G_x^{(p)}$. By
the first part of the proof, this implies that the number of choices for $G_x^{(p)}, G_e ^{(p)}, \alpha_e
|_{G_e^{(p)}}$ and $\alpha_{\bar{e}}|_{G_e^{(p)}}$ is at most $p^{\frac {k_0 ^2 +k_0 +2}{2}}$.

Since all the other $G_e^{(p_i)}$ and $G_x^{(p_i)}$ have fixed cardinality, we have a constant total number of
choices for them and the injections $\alpha_e |G_e^{p_i}$, say $c_0$. Recall that once all the $G_x^{(p_i)}$'s
and $G_e^{(p_i)}$'s are chosen, the number of $G_x$ with a given fixed Sylow system is at most $(pn)^{75
\mu(pn)+16}$(\cite{P}). Recall also that the injections $\alpha_e$ are determined by its restriction to Sylow
subgroups of $G_e$ since they generate the group. Finally we have the following upper bound.
\begin{equation*}
u_n \leq c_0 p^{\frac{(k_0+1)^2+(k_0+1)+2}{2}}(pn)^{75 \mu(pn)+16} \leq c_0
(c_1)^{\frac{k_0^2+5k_0+8}{2}}(pn)^{75\mu+16}
\end{equation*}
where $c_1=p$ and $\mu =\mu(pn)$. \;\;\;\;\qed

\medskip
{\it Remark.} It follows from the proof above that each prime factor of $|G_x|$ is less than or equal to $p$,
thus in the case $p=2$, $u(n)=0$ if $n$ is not a power of 2.

\small{

Yale University,
New Haven, CT 06520-8283
USA

seonhee.lim@yale.edu

and

E.N.S. Paris, 45 rue d'Ulm, 75230 Paris Cedex, France

Seonhee.Lim@ens.fr
\end{document}